\documentclass[10pt]{article}
\usepackage{latexsym}
\usepackage{amsfonts,amsmath,amssymb}
\usepackage{graphics, epsfig}
\usepackage[graphicx]{realboxes}
\usepackage[latin1]{inputenc}
\usepackage{tabularx}
\usepackage{ucs}
\usepackage{marvosym} 
\usepackage[T1]{fontenc}
 \usepackage[english]{babel}
\usepackage[T1]{fontenc}            
\usepackage{textcomp}   
\usepackage{color}
\usepackage{amssymb}
\usepackage{empheq}
\usepackage[colorlinks=true, pdfborder={0 0 0}]{hyperref}
\usepackage{empheq}
\usepackage[latin1]{inputenc}
\usepackage{amsbsy,amsmath,amssymb,amsfonts,latexsym,epsfig,textcomp,wasysym}%
\usepackage{mathtools}
\usepackage{amsmath}
\usepackage{subfigure} 
\usepackage{sidecap}
\usepackage{float} 
\usepackage{marvosym} 
\usepackage[T1]{fontenc} 
\usepackage[english]{babel}
\usepackage[font=small,labelfont=bf,tableposition=top]{caption}
\renewcommand{\thesection}{\arabic{section}}
\renewcommand{\theequation}{\thesection.\arabic{equation}}
\def\R{\mathbb R} \def\N{\mathbb N} 

\usepackage{amsmath,amssymb,amsthm,amsfonts}
\usepackage{enumerate,color}
\makeatletter
\renewcommand{\section}{%
\@startsection{section}{1}{\z@}
{0.5truecm plus -1ex minus -.2ex}%
{1.0ex plus .2ex}{\bfseries\large}}
\def
\@seccntformat#1{\csname the#1\endcsname.\ }
\makeatother
\numberwithin{equation}{section} 
\newtheorem{thm}{Theorem}[section]

\newtheorem{lem}[thm]{Lemma}
\newtheorem{remark}[thm]{Remark}
\newtheorem{proposition}[thm]{Proposition}
\theoremstyle{definition}

\newcommand{\ep}{\varepsilon}
\newcommand{\pa}{\partial}
\newcommand{\Rn}{\mathbb{R}^n}

\newcommand\Label[1]{&\refstepcounter{equation}(\theequation)\ltx@label{#1}&}
\begin{document}
%
%
\begin{NoHyper}
\footnote[0]
	{
	$^{\sharp}$Corresponding author: giuseppe.viglialoro@unica.it
	}
\end{NoHyper}
\begin{center}
	\Large{{\bf 
A refined criterion and lower bounds for the blow--up time in a parabolic--elliptic chemotaxis system with nonlinear diffusion 
   }}
\end{center}
\vspace{5pt}
\begin{center}
	Monica Marras$^1$, 
Teruto Nishino
and Giuseppe Viglialoro$^{1,\sharp}$
\end{center}
\begin{center}
$^1$Dipartimento di Matematica e Informatica
\\  
 Universit\`a di Cagliari, Viale Merello 92, 09123 (Italy)
\end{center}
\vspace{2pt}
\newenvironment{summary}
{\vspace{.5\baselineskip}\begin{list}{}{%
	\setlength{\baselineskip}
	{0.85\baselineskip}
	\setlength{\topsep}{0pt}
	\setlength{\leftmargin}{12mm}
	\setlength{\rightmargin}{12mm}
	\setlength{\listparindent}{0mm}
	\setlength{\itemindent}{\listparindent}
	\setlength{\parsep}{0pt}
	\item\relax}}{\end{list}
	\vspace{.5\baselineskip}}
\begin{abstract}
This paper deals with unbounded solutions to the following zero--flux chemotaxis system
\begin{equation}\label{ProblemAbstract}
\tag{$\Diamond$}
	\begin{cases}
		u_t=\nabla \cdot [(u+\alpha)^{m_1-1} 
		\nabla u-\chi u(u+\alpha)^{m_2-2} 
		\nabla v] 
		&
		(x,t) \in \Omega \times (0,T_{max}), 
		\\[1mm]
		0=\Delta v-M+u 
		& 
		(x,t) \in \Omega \times (0,T_{max}), 
	\end{cases}
\end{equation}
where $\alpha>0$, $\Omega$ is a smooth and bounded domain of $\R^n$, with $n\geq 1$, $t\in (0, T_{max})$, where $T_{max}$ the blow-up time, and $m_1,m_2$ real numbers. Given a sufficiently smooth initial data $u_0:=u(x,0)\geq 0$ and set $M:=\frac{1}{|\Omega|}\int_{\Omega}u_0(x)\,dx$, from the literature it is known that under a proper interplay between the above parameters $m_1,m_2$ and the extra condition $\int_\Omega v(x,t)dx=0$,  system \eqref{ProblemAbstract} possesses for any $\chi>0$  a unique classical solution  which becomes unbounded at $t\nearrow T_{max}$. In this investigation we first show that for  $p_0>\frac{n}{2}(m_2-m_1)$ any blowing up classical solution in $L^\infty(\Omega)$--norm blows up also in $L^{p_0}(\Omega)$--norm. Then we estimate the blow--up time $T_{max}$ providing a lower bound $T$.
\end{abstract}
\section{Introduction and motivation}\label{Section1Introduction}
In this paper we study properties of given solutions which classically solve this chemotaxis problem 
\begin{equation}\label{sys1}
	\begin{cases}
		u_t=\nabla \cdot [(u+\alpha)^{m_1-1} 
		\nabla u-\chi u(u+\alpha)^{m_2-2} 
		\nabla v] 
		&
		(x,t) \in \Omega\times (0,\infty), 
		\\[1mm]
		0=\Delta v-M+u 
		& 
		(x,t) \in \Omega \times (0,\infty), 
		\\[2mm] 
		u_\nu 
		= 
		v_\nu 
		= 0 
		& (x,t) \in \partial\Omega \times 
		(0,\infty), 
		\\[1mm] 
		u(x, 0)=u_0(x)
		& x \in \Omega, 
		\\[1mm] 
\displaystyle		\int_\Omega v(x,t)dx=0 &  t\in (0,\infty), 
	\end{cases}
\end{equation}
where $\alpha, \chi>0$, the spatial variable $x$ is a vector of $\R^n$, with $n\geq 1$, belonging to a smooth and bounded domain $\Omega$ and $t$ is the time variable. Further $m_1,m_2$ are proper real numbers, $\nu$ is the outward normal vector 
to $\pa\Omega$ and the initial data $u_0:=u_0(x)$, supposed to be nonnegative and sufficiently regular, defines also the constant $M$ through the relation $M=\frac{1}{|\Omega|}\int_\Omega u_0(x)dx$. 

In the framework of self organization mechanisms for biological populations, and similarly to many variants of the well--known Keller--Segel models (see the celebrated papers \cite{K-S-1970,Keller-1971-MC,Keller-1971-TBC}), system \eqref{sys1}, which is expressed as a particular case of a more general formulation provided in \cite{WinDj}, represents the situation where the motion of a certain cell density $u(x, t)$ at the position $x$ and at the time $t$, living in an impenetrable (homogeneous Neumann boundary conditions) domain and initially distributed according to the law of $u_0(x)$, is influenced by the presence of a chemical signal concentrations, whose deviation from its spatial mean at the same position time is indicated with $v(x, t)$.  
\begin{remark}
Let us precise that in this paper the mentioned deviation $v$ is, essentially, the difference between the signal concentration and its mean, and that conversely to what happens to the cell and signal densities (which are nonnegative) it changes sign. In particular, from the definition itself of $v$, we have that its mean is zero (as fixed in the last assumption of problem \eqref{sys1}), which in turn ensures the uniqueness of the solution for the Poisson equation under homogeneous Neumann boundary conditions. In this sense, the corresponding compatibility condition, leads for $t>0$ to $\int_\Omega u (x,t)dx=M|\Omega|$ and  by virtue of $\int_\Omega u(x,t)dx=\int_\Omega u_0(x)dx$ (coming by integrating over $\Omega$ the equation for $u$), the choice $M=\frac{1}{|\Omega|}\int_\Omega u_0(x)dx$ remains justified. Finally, and in line with what we said, we advise the reader that in the literature, and also in different places of this section, $v$ stands for the chemical signal concentration itself and not for its deviation; nevertheless, we understand that in view of what we specified above it is not necessary to introduce a further symbol for the deviation, since it will be very clear from the context to which one of these quantities we are referring to.  
\end{remark}
A natural and singular situation possible appearing also in more general cellular processes than that introduced in  model \eqref{sys1}, but also idealized by two partial differential equations  (one for the cell distribution and one for the chemical), is the \textit{chemotaxis collapse}, when an uncontrolled gathering of cells at certain spatial locations is perceived as time evolves;  essentially, $u$, in a particular instant of time (the blow-up time), becomes unbounded in one or more points of its domain. This degeneration of the cell movements into aggregation is, above all, justified by the presence of the destabilizing effect in the coupled term (cross-diffusion term in the evolutive equation for $u$: in our case the expression $\chi u(u+\alpha)^{m_2-2}\nabla v$ in system \eqref{sys1}); in turn the strength of such a destabilizing factor depends on the evolutive equation of $v$ (in our case, of course the second one in system \eqref{sys1}).

The pioneer Keller--Segel system \cite{Keller-1971-MC}, already cited, is obtained from \eqref{sys1} when $m_1=1$, $m_2=2$ and with second equation given by $\tau v_t=\Delta v -v+u$, with $\tau\in\{0,1\}$, where in this case $v$ is the chemical signal concentration (and not its deviation). For positive chemical and cell distributions the expression   $-v+u$ manifests how an increase of the cells favors a production of the signal. For this  case a very comprehensive and extensive theory on existence and properties of global, uniformly bounded or  blow-up solutions, especially in terms of the size of the initial data, is available; for a complete picture, we suggest the introduction of \cite{HorstWink} for the parabolic-parabolic case (i.e., $\tau=1$), \cite{JaLu} and \cite{Nagai} for the parabolic-elliptic case (i.e., $\tau=0$) and in addition the survey by \cite[Hillen and Painter]{Hillen2009UGP} where, inter alia, reviews of various models about Keller-Segel-type systems are discussed.

Besides the size of the initial data, there is another aspect related to the existence of both bounded or unbounded solutions to chemotaxis--systems; this is the mutual interplay between the weight of diffusion $A(u,v)$ and that of the chemotactic sensitivity  $B(u,v)$, which in our context are $(u+\alpha)^{m_1-1}$ and $\chi (u+\alpha)^{m_2-2}$, respectively. (To the readers interested to numerical simulations indicating the influence of the parameters $m_1$ and $m_2$ on solutions to a Keller--Segel system similar to that studied in the current research we suggest \cite[$\S$5]{ViglialoroWoolleyParaelipLogistic}.) Let us give some information in this regard for  system \eqref{sys1} with second equation of parabolic type, i.e  $v_t=\Delta v -v+u$. In \cite{CieslakStinnerFiniteTimeBlowUp}, \cite{CieslakStinnerNewCritical} and  \cite{TaoWinkParaPara} it is essentially established that the relation $m_2<m_1+\frac{2}{n}$ is a necessary and sufficient condition to ensure global existence and boundedness of solutions even emanating from large initial data. This is a generalization of \cite[Theorems 4.1 and 6.1]{HorstWink}, where $m_1=1$ (see also \cite{TaoWinkParaPara} and \cite{ISHIDASekiYokota}). Even more, in \cite{ISHIDAYokota} a parabolic--parabolic degenerate chemotaxis system ($\alpha=0$ in \eqref{sys1}) is discussed: resorting to the natural concept of weak solutions, it is shown that for  $m_2<m_1+\frac{2}{n}$ and $\Omega=\R^n$ the problem possesses global bounded solutions (we refer also to \cite{ISHIDAYokotaSmallData} for a discussion on the super--critical case $m_2\geq m_1+\frac{2}{n}$). Let us note that the reciprocal iteration involving $m_1,m_2$ and $n$ somehow  establishes that the destabilizing effect of the chemo-sensitivity $B(u,v)$ is weaker than that from the diffusion $A(u,v)$, which conversely tends to provide equilibrium to the model. 

Motivated by the above discussion, aim of the present research is expanding the
theory of the mathematical analysis of problem \eqref{sys1} studied in \cite{WinDj}, which, so far we are aware, covers the following situations: (i) for $m_1\leq 1$, $m_2<m_1+\frac{2}{n}$, any sufficiently regular initial data emanates solutions which are global and uniformly bounded; (ii) for $m_1\leq 1$, $m_2>1$, $m_2>m_1+\frac{2}{n}$ and $\Omega$ a ball of $\R^n$ there exist initial data $u_0$ which emanates unbounded solutions at some finite time $T_{max}$. 

In light with this, we are interested in deriving a lower bound $T$ for the blow--up time $T_{max}$ of the unbounded solutions to \eqref{sys1}, so to essentially obtain a safe interval of existence $[0, T)$ where such solutions exist. We will achieve this result according to the steps specified in the next section. 
\section{Some premises and preparatory tools: plan of the paper}\label{SectionPremises}
For these coming reasons, we want to observe that there is no automatic connection between the occurrence of blow-up for solutions to \eqref{sys1} in the classical $L^\infty(\Omega)$--norm and that in $L^p(\Omega)$--norm ($p>1$). Indeed, once it is assumed that $\Omega$ is a bounded domain, we only can conclude that
\[\| u(\cdot,t)\|_{L^p(\Omega)}\leq |\Omega|^\frac{1}{p}\|u(\cdot,t)\|_{L^\infty(\Omega)},\]
so that if a solution blows up in $L^p(\Omega)$--norm, it does in $L^\infty(\Omega)$--norm; conversely, if a solution becomes unbounded in  $L^\infty(\Omega)$--norm at some finite time $T_{max}$, $\int_\Omega u^p$ might also remain bounded in a neighborhood of $T_{max}$. In particular, since for a classical solution $(u,v)$ to system \eqref{sys1} the $u$--component is continuous in $I=[0,T_{max})$, the function $\int_\Omega u^p$  enjoys this same property on $I$ and if $\limsup \int_\Omega u^p$  is finite as $t\searrow T_{max}$, $\int_\Omega u^p$ can even  be continuously prolonged up to the boundary $T_{max}$.

On the other hand,  the evolution in time for the function $t \mapsto \int_\Omega u^p$ is more amenable to be analyzed than that for $t \mapsto \| u(\cdot,t)\|_{L^\infty(\Omega)}$, so that it is preferable to use this to derive lower bounds for $T_{max}$. In this way, in order to avoid the gap between the analysis of the blow--up time $T_{max}$ in the two different mentioned norms will move toward a twofold action: 
\begin{itemize}
\item to detect proper $L^p$--norms, for suitable $p$ depending on $n,m_1$ and $m_2$, ensuring that the unbounded solutions in $L^\infty(\Omega)$--norm also blow up in these $L^p(\Omega)$--norms; 
\item to provide lower bounds for the blow--up time of unbounded solutions in these $L^p(\Omega)$--norms.
\end{itemize}
To be more precise, we invoke the results in \cite{WinDj} in order to frame scenarios where local classical solutions $(u,v)$ to system \eqref{sys1} are detected ($\S$\ref{SectionStartingAndMainTheores}). Successively, and this is a crucial step in our investigation, we show that under suitable assumptions on the parameters $m_1$ and $m_2$ any local solution to system \eqref{sys1} which blows up at finite time $T_{max}$ in $L^\infty(\Omega)$--norm also does in $L^p(\Omega)$--norm (Theorem \ref{TheoremMain1} of $\S$\ref{SectionStartingAndMainTheores}, proved in $\S$\ref{SectionEnergyFunction}); to this aim, we will rely on \cite[Theorem 2.2]{FREITAGLpLinftyCoincide}, so to derive proper estimates by virtue of the analysis of the energy function $\Phi(t)=\frac{1}{p}\int_\Omega (u+\alpha)^p$, for some $p>1$, defined for all $t\in (0,T_{max})$ and associated to the local solution $(u,v)$; these estimates are derived in  $\S$\ref{SectionEnergyFunction}. (In $\S$\ref{PreliminariesSection} we give some necessary and preliminary tools.) Successively, in $\S$\ref{SectionODIForLower}, it is established that the same  $\Phi(t)$ satisfies a first order differential inequality
(ODI) of the type $\Phi'(t)\leq \Psi(\Phi(t))$ on $(0,T_{max})$. In particular, for any $\tau>0$ the function $\Psi(\tau)$
obeys the Osgood criterion (\cite{Osgood}),
\begin{equation}\label{OsgoodCriterion}
\int_{\tau_0}^\infty \frac{d\tau}{\Psi(\tau)}<\infty\quad \textrm{with } \tau_0>0,
\end{equation}
so that an integration on $(0,T_{max})$ of the mentioned ODI 
implies $$T_{max}\geq \int_{\Phi(0)}^\infty \frac{d \Phi}{\Psi(\Phi)}:=T,$$ thereby yielding the desired lower bound $T$ for the blow--up time $T_{max}$. (This is Theorem \ref{main} of $\S$\ref{SectionStartingAndMainTheores}, whose proof is presented in $\S$\ref{SectionODIForLower}.)
\section{Starting point and presentation of the main theorems}\label{SectionStartingAndMainTheores}
From the above considerations, let us give the following proposition, which represents the starting point of our work and that we claim according to our purposes. 

First,  we fix these mutual \textit{blow--up restrictions} on the parameters $m_1,m_2$, since in the light of the results presented in $\S$\ref{Section1Introduction} they are the natural assumptions enforcing solutions to model \eqref{sys1} to become unbounded:
\begin{align}
\label{b}\tag{\textsf{BU}}
m_2>m_1+\frac{2}{n},  \quad m_1\le 1,\quad  m_2 > 1.
\end{align}
\begin{proposition}\label{WinklerDjiePropo}
Let $\Omega$ be a bounded and smooth domain of $\R^n$, with $n\geq 1$, $\alpha,\chi>0$ and $0\leq u_0\in C^\kappa(\bar{\Omega})$, for some $\kappa>0$, a nontrivial initial data with $M=\frac{1}{|\Omega|}\int_\Omega u_0(x)dx$. Additionally, let $m_1,m_2 \in \R$ comply with the \textit{blow--up restrictions} \eqref{b}. Then, there exist a finite time $T_{max}>0$ and a unique local classical solution 
\begin{align*}
	&
	(u,v) \in C\left(\bar{\Omega} 
	\times [0,T_{max})\right)
	\cap 
	C^{2,1}\left(\bar{\Omega} 
	\times (0,T_{max})\right)\times
	C^{2,0}\left(\bar\Omega
	\times (0,T_{max})\right)
\end{align*}
to system \eqref{sys1} which blows up at  $T_{max}$ in the sense that 
\begin{equation}\label{Blow-UpClassicalSense}
\limsup_{ t \nearrow  T_{max }}\lVert u(\cdot, t)\rVert
_{L^{\infty}(\Omega)} =
\infty. 
\end{equation}
\begin{proof}
See \cite[Theorem 4.5]{WinDj}.
\end{proof}
\end{proposition}
\begin{remark}
For the sake of scientific information, \cite[Theorem 4.5]{WinDj} is proved in a ball of $\R^n$ and moreover under additional restrictions on the data $u_0$, as in particular some  assumptions on its support which rule out the choice of constant initial data. (Indeed, $(u,v)=(constant,0)$ is a bounded global solutions to system \eqref{sys1}; this is the reason why we exclude trivial $u_0$ in Proposition \ref{WinklerDjiePropo} and throughout all the paper.) Despite that, since in the present investigation we are mostly interested in the derivation of lower bounds for the blow--up time $T_{max}$ to unbounded solutions to system \eqref{sys1}, we understand that the more general claim proposed in Proposition \ref{WinklerDjiePropo} does not mislead and is consistent with our overall aim.   
\end{remark}
\begin{thm}\label{TheoremMain1} 
Let $\Omega$ be a bounded and smooth domain of $\R^n$, with $n\geq 1$,  $\alpha,\chi>0$ and $0\leq u_0\in C^\kappa(\bar{\Omega})$, for some $\kappa>0$,  a nontrivial initial data. Then, for $m_1,m_2 \in \R$ complying with the \textit{blow--up restrictions} \eqref{b} and $M=\frac{1}{|\Omega|}\int_\Omega u_0(x)dx$,
the blow--up classical solution $(u,v)$ to system \eqref{sys1} provided by Proposition \ref{WinklerDjiePropo} is such that for all $p_0>\frac{n}{2}\left(m_2-m_1\right)$
\begin{align*}
\limsup_{t \nearrow T_{\rm max}}\, 
\left\|u(\cdot,t)\right\|_{L^{p_0}(\Omega)} 
= \infty.
\end{align*}
\end{thm} 
\begin{thm}\label{main}
Let $\Omega$ be a bounded and smooth domain of $\R^n$, with $n\geq 1$,  $\alpha,\chi>0$ and $0\leq u_0\in C^\kappa(\bar{\Omega})$, for some $\kappa>0$,  a nontrivial initial data. Then, for $m_1,m_2 \in \R$ complying with the \textit{blow--up restrictions} \eqref{b} and $M=\frac{1}{|\Omega|}\int_\Omega u_0(x)dx$, it is possible to find  $\bar{p}>1$ and  $E_5, E_8, E_9>0$ as well as $\gamma,\delta>1$, depending on $\bar{p}$, such that the blow--up time $T_{max}$ of the unbounded classical solution $(u,v)$ to system \eqref{sys1} provided by Proposition \ref{WinklerDjiePropo} satisfies
\begin{align}\label{blow-up-time}
		T_{max} \geq \int_{\Phi(0)}^{\infty}  \frac{d\tau}  {E_8{\tau}^{\gamma} + E_9{\tau}^{\delta} + E_5},
\end{align}
where $\Phi(0)= \frac{1}{\bar{p}}  \int_{\Omega} \left(u_0+  \alpha\right)^{\bar{p}}$.
\end{thm}
\begin{remark}\label{RemarkOnSignOfLowerBound}
In the absence of the result of Theorem \ref{TheoremMain1},  and taking in mind what discussed at the beginning of $\S$\ref{SectionPremises}, the current formulation of Theorem \ref{main} might fail without adding the extra hypothesis that $\limsup_{t\rightarrow T_{max}}\Phi(t)=\infty.$ In fact, the above ODI $\Phi'(t)\leq \Psi(\Phi(t))$ would infer, by integration on $(0,T_{max})$ as well, that 
\[T_{max}\geq \int_{\Phi(0)}^{{\Phi(T_{max})}} \frac{d\tau}{\Psi(\tau)},\]
which does not produce any lower bound if no additional assumption on $\Phi(T_{max})$ is given. Thereafter, even though in the literature there are several papers concerning estimates for lower bounds of blow-up time for solutions to general evolutive problems whose formulation relies on the hypothesis on the divergence of certain energy functions  (see, for instance, \cite[Theorem 1 and Theorem 2]{LI-ZhengBlowUpK-S},  \cite[Theorem 2.4 and Theorem 2.7]{MarrasVernierVigliaWithm} and \cite[Theorem 1 and Theorem 2]{PSong} for contributions in the frame of chemotaxis models or \cite[Theorem 2.1]{PAYNEPhilSchaferNonlinAnalysis} and \cite[Theorem 1 and Theorem 4]{PAYNE-Phil-SchaferBounds} for others in different areas), the inspiring paper \cite{FREITAGLpLinftyCoincide} represents a cornerstone that allows us  to avoid this hypothesis, precisely thanks to the implication 
\begin{center}
``A solution to \eqref{sys1} which blows up in $L^\infty(\Omega)$--norm automatically does in $L^p(\Omega)$--norm''
\end{center}
given in Theorem \ref{TheoremMain1}. (An equivalent approach is employed in \cite{NishinoYokotaclassicalCoincideLp} for unbounded solutions to the same fully parabolic chemotaxis problem analyzed in \cite{FREITAGLpLinftyCoincide}.)
\end{remark}
\section{Fixing some parameters and functional inequalities}\label{PreliminariesSection}
In the following lemma, we fix the value of an important parameter, used to quantify certain constants appearing throughout our logical steps, essentially by adjusting the data $m_1,m_2$ and $n$ defining problem \eqref{sys1}. This parameter will be set in a such a way that the employments of some crucial inequalities below will be straightforwardly justified. 
\begin{lem}\label{Lemmapbarra}  
For $n\in \N$, let $m_1, m_2$ satisfy  the assumptions in \eqref{b} and $p_0>\frac{n}{2}(m_2-m_1)$. Additionally, for any $q_1>n+2$, $q_2>(n+2)/2$, let
\begin{equation}\label{ConstantForTechincalInequality_Barp}
\bar{p}:=\max 
\begin{Bmatrix}
p_0 \vspace{0.1cm}   \\ 1-m_1 \vspace{0.1cm} 
 \\ 3-m_2 \vspace{0.1cm} 
 \\ 2-m_1-\frac{2}{n} \vspace{0.1cm} \\
p_0-m_2+1 \vspace{0.1cm} \\
q_1\vspace{0.1cm}
 \\q_1(m_2-1) \vspace{0.1cm}\\ 1-m_1\frac{(n+1)q_1-(n+2)}{q_1-(n+2)} \vspace{0.1cm}  \\ 1-\frac{m_1}{1-\frac{n}{n+2}\frac{q_2}{q_2-1}} 
\end{Bmatrix}
+1.
\end{equation}
Then for all $p \geq  \bar{p}$ these relations hold
\[
\begin{subequations}
  \begin{tabularx}{\textwidth}{Xp{1.8cm}X}
 \begin{equation}
    \label{2}
p> \frac n2 (1-m_1)
  \end{equation}
  & &
  \begin{equation}
    \label{InequalityForG-NInu^p+1}
0< a_1:=\frac{\frac{(p+m_1-1)}{2}(1-\frac{1}{p})}{\frac{(p+m_1-1)}{2}+\frac{1}{n}-\frac{1}{2}}<1
  \end{equation}
  \\
 \begin{equation}
    \label{CoeffBeta}
0 <\beta:=\frac
{\frac{p+m_2-1}{2p_0}
-\frac{1}{2}}
{\frac{p+m_1-1}{2p_0}+\frac{1}{n}-
\frac{1}{2}}<1
  \end{equation}
  & &
    \begin{equation}\label{5}
\frac{1}{k}>\frac{1}{2}-\frac{1}{n}
\end{equation}
  \\
  \begin{equation}
    \label{Coeff_b1}
0 <a_2:=\frac
{\frac{p+m_1-1}{2p_0}
-\frac{1}{k}}
{\frac{p+m_1-1}{2p_0}+\frac{1}{n}-
\frac{1}{2}} <1
  \end{equation}
& & 
\begin{equation}\label{a_1}
0<a_3:=\frac{\frac{p+m_1-1}{2p}-\frac{1}{k}} {\frac{p+m_1-1}{2p}+\frac{1}{n}-\frac{1}{2}}<1
\end{equation}
\\
\begin{equation}\label{rho<1andgamma>delta>1}
0<\sigma<1\quad \textrm{ and } \quad \gamma > \delta >1,
\end{equation}
&&
  \end{tabularx}
\end{subequations}
\]
where
\begin{equation}\label{gammabetadelta}
k:=\frac{2(p+m_2-1)}{p+m_1-1},\quad \delta= \frac{p+m_2-1}{p},\quad \sigma= \frac{k a_3}{2},\quad \gamma = \frac{p+m_2-1}{p}\frac{1-a_3}{1-\sigma}.
\end{equation}
\begin{proof}
The assumptions done in \eqref{b} and the definition of $\bar{p}$, in conjunction with the restriction on $p_0$, imply $\bar{p}>\frac{n}{2}(1-m_1)$, i.e. \eqref{2}, so that in turn  we have
\begin{equation*}
\frac{p+m_1-1}{2}\geq \frac{\bar{p}+m_1-1}{2}> \frac{\bar{p}+m_1-1}{2\bar{p}}\quad \textrm{ and that } \quad \frac{\bar{p}+m_1-1}{2\bar{p}}>\frac{n-2}{2n};
\end{equation*}
therefore
\begin{equation*}
1-\frac{n}{2}+\frac{n}{2}(p+m_1-1)>0
\end{equation*}
and, thus, also \eqref{InequalityForG-NInu^p+1} is attained. Again from the assumption on $p_0$, we also have, recalling again \eqref{b}, that $\left(1-\frac{2}{n}\right)p_0+1-m_1>0$, so that relation $p>{p_0}-m_2+1$ and the definition of $k$ also easily give  \eqref{CoeffBeta} and \eqref{5}, this last one also used to show  \eqref{Coeff_b1} and \eqref{a_1}. The remaining inequalities come from 
$$\frac{p+m_1-1}{2p}+\frac{1}{n}-
\frac{1}{2}>0, \quad k>2 \quad \textrm{and}\quad 1-\sigma=1-\frac{k a_3}{2}<(1-a_3).$$
\end{proof}
\end{lem}
\begin{remark}\label{RemarkPbar}
As to the definition of the parameter $\bar{p}$ in \eqref{ConstantForTechincalInequality_Barp}, we desire to point out that  the expression of $\bar{p}$ is precisely fixed in that way exactly to avoid to have to enlarge a general $p>1$ up to some suitable values which are used in our derivations. Moreover, the addition ``+1'' in the same definition is not strictly necessary but, undoubtedly, its presence will allow us to uniquely establish the magnitudes of some constants, many of these taking part, inter alia, in the quantitative calculations of the lower bound $T$ for $T_{max}$ of Theorem \ref{main}.
\end{remark}
%
As announced, let us now recall the Gagliardo--Nirenberg inequality, which throughout this paper will be used in a less common version:
\begin{lem}\label{G-N_InequLemma}
Let $\Omega$ be a  bounded and smooth domain of $\Rn$, with $n\geq 1$, and $m_1,m_2\in \R$ complying with the \textit{blow--up restrictions} \eqref{b}. Additionally, for $p=\bar{p}$ given in \eqref{ConstantForTechincalInequality_Barp}, let 
	 $\mathfrak{q},\mathfrak{s}\in [\frac{2}{p+m_1-1},\frac{2p}{p+m_1-1}]$, $\mathfrak{p}\in [\frac{2p}{p+m_1-1},\frac{2(p+m_2-1)}{p+m_1-1}]$. Then there exists a uniquely determined positive constant  
	$C_{GN}=C_{GN}(n,m_1,m_2,\Omega)$ such that 
	\begin{align*}
		\|w\|_{L^\mathfrak{p}(\Omega)}\leq C_{GN}
		\left(\|\nabla w\|_{L^2(\Omega)}
		^{a}\|w\|_{L^\mathfrak{q}(\Omega)}^{1-a}
	+\|w\|_{L^\mathfrak{s}(\Omega)}\right)\quad \textrm{ for all } w\in W^{1,2}(\Omega)\cap 
	L^\mathfrak{q}(\Omega),
	\end{align*}
	where 
$a:=\frac{\frac{1}{\mathfrak{q}}-\frac{1}{\mathfrak{p}}}
		{\frac{1}{\mathfrak{q}}+\frac{1}{n}-
		\frac{1}{2}}\in (0,1)$. 
\begin{proof}
This is an adaptation and a specific case of \cite[Lemma 2.3]{Li-Lankeit_2016} with $\mathfrak{r}=2$. Given for $\bar{p}$ as in \eqref{ConstantForTechincalInequality_Barp}, for $\mathfrak{p}$ and $\mathfrak{q}$ as in our assumptions, we see that  $0<\mathfrak{q} \leq \mathfrak{p}\leq\infty$ and, by virtue of \eqref{2},  $\frac{1}{\mathfrak{r}}\leq \frac{1}{n}+\frac{1}{\mathfrak{p}}$, so that the claim straightforwardly follows by  \cite[Lemma 2.3]{Li-Lankeit_2016}, taking as $C_{GN}$ the maximum value of the constant $c$ therein used under the constraints on $\mathfrak{p}, \mathfrak{q}, \mathfrak{r}$ and $\mathfrak{s}$ herein assumed. 
\end{proof}
\end{lem}
\section{The energy function $\Phi(t) := \frac{1}{p}\int_\Omega (u+\alpha)^{p}$: some a priori estimates}\label{SectionEnergyFunction}
Having ensured existence of unbounded solutions $(u,v)$ to system \eqref{sys1}  we can now turn our attention to the evolution in time of the energy function $\Phi(t) := \frac{1}{p}\int_\Omega (u+\alpha)^{p}$, with $p>1$. (In particular, having properly fixed the parameter $\bar{p}$ in Lemma \ref{Lemmapbarra}, and in light of Remark \ref{RemarkPbar}, we will analyze $\Phi(t)$ for $p=\bar{p}.$) Apparently, this analysis will provide crucial
information for both the proofs of Theorems \ref{TheoremMain1} and \ref{main}. 
\begin{lem}\label{LemmaFirstStepDerivativePhi}
Under the assumptions of Proposition \ref{WinklerDjiePropo}, let $(u,v)$ be the local solution to system \eqref{sys1} which blows up at finite time $T_{max}$ in the sense of \eqref{Blow-UpClassicalSense}, and $\Phi(t)$ the energy function defined for $p=\bar{p}$ as in \eqref{ConstantForTechincalInequality_Barp} by
\begin{equation}\nonumber
 \Phi(t) := \frac{1}{p} 
		\int_{\Omega} \left(u+\alpha\right)^p  \quad \textrm{ on } (0,T_{max}).
\end{equation}
Then there exist $E_0,E_1,E_5>0$ such that 
\begin{align}\label{ph}
\Phi'(t) 
&\le  - E_0 \int_{\Omega} \left| \nabla (u+\alpha)^{\frac{p+m_1-1}{2}}\right|^2 
+E_1\int_{\Omega}(u+\alpha)^{p+m_2-1}
+E_5\quad \textrm{for all } t\in (0,T_{max}). 
\end{align} 
\begin{proof}
The first equation of \eqref{sys1} enables us to see 
\begin{align}\nonumber
		\frac{1}{p}\dfrac{d}{dt}\int_{\Omega}(u+
\alpha)^p
		&=-\int_{\Omega}\nabla (u+\alpha)^{p-1} 
\cdot
		\left[(u+\alpha)^{m_1-1}\nabla u
		-\chi u(u+\alpha)^{m_2-2}
		\nabla v
\right]
		\\\label{FirstDerivativePhiPrevious}
&=:-I_1 + I_2\quad \textrm{ for all } t\in(0,T_{max}), 
	\end{align}
where 
\begin{align*}
&I_1:=\int_{\Omega}\nabla (u+\alpha)^{p-1} 
\cdot (u+\alpha)^{m_1-1}\nabla u\quad \textrm{ for all } t\in (0,T_{max}), 
\\
&I_2:=\int_{\Omega}\nabla (u+\alpha)^{p-1}
\cdot \chi u(u+\alpha)^{m_2-2} 
\nabla v\quad \textrm{ for all } t\in (0,T_{max}). 
\end{align*}
As to the addendum $I_1$, from Lemma  \ref{Lemmapbarra} we have $p>1-m_1 $, so that we can write 
\begin{align}\nonumber
 I_1&= \int_{\Omega}\nabla (u+\alpha)^{p-1} 
\cdot (u+\alpha)^{m_1-1}\nabla u=(p-1)\int_{ \Omega} 
(u+\alpha)^{p+m_1-3}|\nabla u|^2\\ \label{I1}
&
=E_0 \int_{ 
 \Omega} \left| \nabla (u+\alpha)^
{\frac{p+m_1-1}{2}}\right|^2\quad \textrm{ for all } t\in (0,T_{max}),
\end{align}
where $E_0=\frac  {4(p-1)} {(p+m_1-1)^2 }.$ Similarly, in order to control $I_2$, we start to write 
\begin{align*}
I_2&=\int_{\Omega}\nabla (u+\alpha)^{p-1}
\cdot \chi u(u+\alpha)^{m_2-2} 
\nabla v=\chi (p-1)\int_{ \Omega} 
u(u+\alpha)^{p+m_2-4}\nabla u \cdot \nabla v  \quad \textrm{ on } (0,T_{max}).  
\end{align*}
Setting 
\begin{align*}
F(u):= \int_{0}^{u} \tau(\tau + \alpha)^{p+m_2-4} d\tau,
\end{align*}
we explicitly see from problem \eqref{sys1} that 
\begin{align}\nonumber
I_2&=\chi (p-1)\int_{ \Omega} 
\nabla F(u) \cdot \nabla v =-\chi (p-1)\int_{\Omega}
F(u) \Delta v=
-\chi (p-1)\int_{\Omega}
F(u) (M-u)\\\label{I2}
&= 
-\chi (p-1)M\int_{\Omega}
F(u) 
+ 
\chi (p-1)\int_{\Omega}
F(u)u 
 \quad \textrm{ on } (0,T_{max}). \
\end{align}
In view again of Lemma \ref{Lemmapbarra}, we have $p>3-m_2$ so 
we can calculate $F(u)$ as 
\begin{align*}
F(u)&= \int_{0}^{u} 
\tau(\tau + \alpha)^{p+m_2-4} d\tau
= 
\frac{1}{p+m_2-3}\Big(\ 
u(u + \alpha)^{p+m_2-3}
-
\int_{0}^{u}(\tau+\alpha)^{p+m_2-3} d\tau\Big). 
\end{align*}
Additionally, from the relation
\begin{align*}
\int_{0}^{u}(\tau+\alpha)^{p+m_2-3} d\tau 
=
\frac{1}{p+m_2-2} 
(u+\alpha)^{p+m_2-2}- 
\frac{1}{p+m_2-2}
\alpha^{p+m_2-2}, 
\end{align*}
with some manipulations and using 
\begin{align}\nonumber 
u(u+\alpha)^{p+m_2-3} 
&=(u+\alpha)^{p+m_2-2}
-\alpha(u+\alpha)^{p+m_2-3},
\end{align}
we infer
\begin{align}\label{FofUExpression}
F(u)
= 
&
\frac
{(u+\alpha)^{p+m_2-2}}{p+m_2-2}
-
\frac
{\alpha(u+\alpha)^{p+m_2-3}}
{p+m_2-3}
+
\frac{\alpha^{p+m_2-2}}
{(p+m_2-3)(p+m_2-2)}. 
\end{align}
Henceforth,  \eqref{I2}-\eqref{FofUExpression} and 
\begin{align}\nonumber 
u(u+\alpha)^{p+m_2-2} 
&=(u+\alpha)^{p+m_2-1}
-\alpha(u+\alpha)^{p+m_2-2}
\end{align}
now produce on $(0,T_{max})$
\begin{align}
\nonumber
I_2 
\nonumber
&=
C_2C_3\int_{\Omega}u(u+\alpha)^{p+m_2-2}
-C_2C_4
\int_{\Omega}u(u+\alpha)^{p+m_2-3}
-C_1C_3\int_{\Omega} 
(u+\alpha)^{p+m_2-2}
\\
\nonumber
&\quad
+C_1C_4\int_{ \Omega} (u+\alpha)
^{p+m_2-3}-C_1C_5|\Omega|+C_2C_5M|\Omega|
\\
%
&=
C_2C_3\int_{\Omega}
(u+\alpha)^{p+m_2-1}
-(\alpha C_2C_3+C_2C_4+C_1C_3)
\int_{\Omega}
(u+\alpha)^{p+m_2-2}
\\\nonumber
&\quad\,
+(\alpha C_2C_4+C_1C_4)\int_{\Omega} 
(u+\alpha)^{p+m_2-3}
-C_1C_5|\Omega|+C_2C_5M|\Omega|
\\
\label{I_2}
&=
E_1\int_{\Omega}(u+\alpha)^
{p+m_2-1}
-E_2
\int_{\Omega}
(u+\alpha)^
{p+m_2-2}
+E_3
\int_{\Omega}
(u+\alpha)^
{p+m_2-3}
+
E_4, 
\end{align}
where we have set
\begin{equation*}
\begin{cases}
C_1=\chi(p-1)M,\quad C_2=\chi(p-1),\\
C_3=\frac{1}{p+m_2-2},\quad
C_4=\frac{\alpha}{p+m_2-3},\quad
C_5=\frac{\alpha^{p+m_2-2}}
{(p+m_2-3)(p+m_2-2)},\\
E_1=C_2C_3,\quad
E_2=\alpha C_2C_3+
C_2C_4+C_1C_3,\\ 
E_3=\alpha C_2C_4 
+C_1C_4,\quad
E_4=-C_1C_5|\Omega|+C_2C_5M|\Omega|. 
\end{cases}
\end{equation*}
On the other hand, since a combination of relations \eqref{I1} and \eqref{I_2} yields 
the following identity
\begin{align}\nonumber
\Phi'(t) 
&=
- E_0
\int_{ 
 \Omega} \left| \nabla (u+\alpha)^
{\frac{p+m_1-1}{2}}\right|^2 
+E_1\int_{\Omega}(u+\alpha)^
{p+m_2-1}
-E_2
\int_{\Omega}
(u+\alpha)^
{p+m_2-2}
\\\label{phb}
&\quad\, 
+E_3
\int_{\Omega}
(u+\alpha)^
{p+m_2-3}
+
E_4\quad \textrm{ for all } t\in (0,T_{max}), 
\end{align}
we can estimate the forth term on the right--hand side of \eqref{phb} by the Young inequality, so to have for any $\delta_0>0$
\begin{align}\label{forth}
E_3 \int_{ \Omega} (u+\alpha)^{p+m_2-3}
\le 
\delta_0 
\int_{ \Omega} (u+\alpha)^{p+m_2-2}
+
D_0(\delta_0)\quad \textrm{ on } (0,T_{max}),
\end{align}
with
\[D_0(\delta_0)=\frac{1}{p+m_2-2}\Big(\delta_0E_3^{-\frac{p+m_2-2}{p+m_2-3}}\frac{p+m_2-2}{p+m_2-3}\Big)^{-(p+m_2-3)}|\Omega|.\]  
Subsequently, from \eqref{phb} and \eqref{forth} is achieved that 
\begin{align}\nonumber
\Phi'(t) 
&\le
- E_0
\int_{ 
 \Omega} \left| \nabla (u+\alpha)^
{\frac{p+m_1-1}{2}}\right|^2 
+E_1\int_{\Omega}(u+\alpha)^
{p+m_2-1}
-(E_2-\delta_0)
\int_{\Omega}
(u+\alpha)^
{p+m_2-2}\\\nonumber
&\quad\, +E_4+D_0(\delta_0) \quad \textrm{ on } (0,T_{max}),
\end{align} 
so that, by taking $\delta_0=E_2$ and considering that from Lemma \ref{LemmaFirstStepDerivativePhi} the constant $E_4$ might be negative,  we finally conclude posing $E_5:=|E_4|+D_0(E_0)$ and obtaining 
\begin{align}\nonumber
\Phi'(t) 
&\le 
- E_0
\int_{ 
 \Omega} \left| \nabla (u+\alpha)^
{\frac{p+m_1-1}{2}}\right|^2 
+E_1\int_{\Omega}(u+\alpha)^
{p+m_2-1}
+E_5 \quad \textrm{ for all } t\in (0,T_{max}).
\end{align} 
\end{proof}
\end{lem}
%
The coming lemma includes the details used to control the terms $\int_{\Omega}|\nabla  (u+\alpha)^\frac{p+m_1-1}{2}|^2$ and  $\int_{\Omega} (u+\alpha)^{p+m_2-1}$. 
\begin{lem}\label{LemmaEstimateu^pandu^p+m2-2} 
Under the assumptions of Proposition \ref{WinklerDjiePropo}, let $(u,v)$ be the local solution to system \eqref{sys1}, which blows up at finite time $T_{max}$ in the sense of \eqref{Blow-UpClassicalSense}, and $\Phi(t)$ the energy function defined in Lemma \ref{LemmaFirstStepDerivativePhi}. Then there exist  positive constants $E_6$ and $\lambda$ such that
\begin{equation}
-\int_{\Omega}
\left|
\nabla (u+\alpha)^{\frac{p+m_1-1}{2}}
\right|^{2}
\le -E_6\Phi^\lambda(t) + 1 \quad \textrm{ for all } t\in (0,T_{max}).
\end{equation}
If, additionally, for some $p_0>\frac{n}{2}(m_2-m_1)$ it is known that 
\begin{equation}\label{AssumptionBoundednessLp0}
\lVert u(\cdot,t)\rVert_{L^{p_0}(\Omega)}\leq L\quad \textrm{ for all } t \in (0,T_{max}),
\end{equation}
then we can find $E_7>0$ such that for all $\varepsilon>0$ it holds
\begin{align}\label{powerc}
\int_{ \Omega} (u+\alpha)^{p+m_2-1}
\le 
\ep
\int_{\Omega}
\left|\nabla(u+\alpha)^{\frac{p+m_1-1}{2}}
\right|^2 
+E_7 \quad \textrm{ on } (0,T_{max}).
\end{align}
\begin{proof}
Thanks to Lemma \ref{Lemmapbarra} we can set 
\[\mathfrak{q}=\mathfrak{s}=\frac{2}{p+m_1-1}, \mathfrak{p}=  \frac{2p}{p+m_1-1}, \] 
and make use of the Gagliardo--Nirenberg 
inequality in Lemma \ref{G-N_InequLemma}. We achieve 
\begin{align}\nonumber
\int_{\Omega} (u+\alpha)^p 
&=
\left\|
(u+\alpha)^{\frac{p+m_1-1}{2}}
\right\|_
{L^{\frac{2p}{p+m_1-1}}(\Omega)}
^{\frac{2p}{p+m_1-1}}
\\\nonumber
&
\le
C_{GN}
\left\|
\nabla (u+\alpha)^{\frac{p+m_1-1}{2}}
\right\|_{L^2(\Omega)}^
{^{\frac{2pa_1}{p+m_1-1}}}
\left\|
(u+\alpha)^{\frac{p+m_1-1}{2}}
\right\|_{L^{\frac{2}{p+m_1-1}}(\Omega)}^
{\frac{2p(1-a_1)}{p+m_1-1}}
\\\nonumber
&
\quad+
C_{GN}
\left\|
(u+\alpha)^{\frac{p+m_1-1}{2}}
\right\|_{L^{\frac{2}{p+m_1-1}}(\Omega)}^
{\frac{2p}{p+m_1-1}}
\\\nonumber
&\le
c_1\left(
\left\|
\nabla (u+\alpha)^{\frac{p+m_1-1}{2}}
\right\|_{L^2(\Omega)}^2
+
1\right)^
{^{\frac{pa_1}{p+m_1-1}}}
\\\label{powera}
&\le
c_1\left(
\int_{\Omega}
\left|
\nabla (u+\alpha)^{\frac{p+m_1-1}{2}}
\right|
1\right)^
{^{\frac{1}{\lambda}}}\quad \textrm{ on } (0,T_{max}),
\end{align}
for some $c_1>0$, ${\lambda}:=\frac{p+m_1-1}{pa_1}>0
$, and where 
\begin{align*}
a_1:=\frac{\frac{p+m_1-1}{2}
(1-\frac{1}{p})}{\frac{p+m_1-1}{2}
+\frac{1}{n}-\frac{1}{2}}
\end{align*}
belongs to $(0,1)$ in view of relation \eqref{InequalityForG-NInu^p+1}. 
%
Subsequently for $E_6:=(\frac{p}{c_1}
)^\lambda$ and through the definition $\int_{\Omega} (u+\alpha)^{p} = p\Phi(t)$, we have that this relation is satisfied
\begin{align}\label{powerb}
-\int_{\Omega}
\left|
\nabla (u+\alpha)^{\frac{p+m_1-1}{2}}
\right|^{2}
&\le 
-
\left(\frac{1}{c_1}
\int_{\Omega} (u+\alpha)^p
\right)^{\lambda}
+1
=-E_6\Phi^\lambda(t) + 1\quad \textrm{ on } (0,T_{max}),
\end{align}
so that the first part of this lemma is shown.

As to the second claim, we will proceed in a similar way to deal with the term $\int_{\Omega} (u+\alpha)^{p+m_2-1}$.  With the aid of bound \eqref{AssumptionBoundednessLp0}, for $k:=\frac{2(p+m_2-1)}{p+m_1-1}$ and Lemma \ref{Lemmapbarra},  if we set   
\[\mathfrak{q}=\frac{2p_0}{p+m_1-1}, \mathfrak{p}= k, \mathfrak{s}=\frac{2}{p+m_1-1}, \] 
the Gagliardo--Nirenberg inequality given in Lemma \ref{G-N_InequLemma} yields  constants $c_3>0$, 
\begin{align*}
a_2&:=
\frac
{\frac{p+m_1-1}{2p_0}
-\frac{1}{k}}
{\frac{p+m_1-1}{2p_0}+\frac{1}{n}-
\frac{1}{2}} 
\in (0,1) \quad (\textrm{recall } \eqref{Coeff_b1}), 
\end{align*} 
and 
\begin{align}\label{point}
\beta&=\frac{ka_2}{2}=\frac
{\frac{p+m_2-1}{2p_0}
-\frac{1}{2}}
{\frac{p+m_1-1}{2p_0}+\frac{1}{n}-
\frac{1}{2}}\in (0,1)\quad (\textrm{recall } \eqref{CoeffBeta}), 
\end{align}
with the property that 
\begin{align}\nonumber
\int_{ \Omega} (u+\alpha)^{{p+m_2-1}}
&=
\left\|
(u+\alpha)^{\frac{p+m_1-1}{2}}
\right\|
_{L^{k}(\Omega)}
^{k}
\\\nonumber
&\le
C_{GN}
\left\|
\nabla (u+\alpha)^{\frac{p+m_1-1}{2}}
\right\|_{L^2(\Omega)}^
{ka_2}
\left\|
(u+\alpha)^{\frac{p+m_1-1}{2}}
\right\|_
{L^{\frac{2p_0}
{p+m_1-1}}(\Omega)}^{
k(1-a_2)}
\\\nonumber
&\quad\, 
+C_{GN}
\left\|
(u+\alpha)^{\frac{p+m_1-1}{2}}
\right\|_{L^{\frac{2}{p+m_1-1}}(\Omega)}^
{\frac
{2(p+m_2-1)}{p
+m_1-1}}
\\\label{power}
&\le
c_3
\left[
1+\left(
\int_{\Omega}
\left|\nabla(u+\alpha)^{\frac{p+m_1-1}{2}}
\right|^2
\right)^\beta\right] \quad \textrm{ on } (0,T_{max}).
\end{align} 
Applying to the gradient term appearing in \eqref{power} the Young inequality, supported with the introduction of an arbitrary positive constant $\varepsilon$, we can write 
\begin{align}\label{Young}
c_3
\left(
\int_{\Omega}
\left|\nabla(u+\alpha)^{\frac{p+m_1-1}{2}}
\right|^2
\right)^\beta
\le 
\ep
\int_{\Omega}
\left|\nabla(u+\alpha)^{\frac{p+m_1-1}{2}}
\right|^2 
+ D_1(\ep)\quad \textrm{ on } (0,T_{max}),
\end{align}
with some $D_1(\ep)>0$,  so as a consequence bound \eqref{power} is reduced to 
\begin{align}\label{powerc}
\int_{ \Omega} (u+\alpha)^{p+m_2-1}
\le 
\ep
\int_{\Omega}
\left|\nabla(u+\alpha)^{\frac{p+m_1-1}{2}}
\right|^2 
+E_7 \quad \textrm{ for all } t\in (0,T_{max}),
\end{align}
where $E_7:=c_3+D_1(\ep)$. 
\end{proof}
\end{lem}
This following result will be the last step toward the proof of Theorem \ref{TheoremMain1}.
\begin{lem}\label{LemmaBoundednessNormUpCrossDiffusion} 
Under the assumptions of Proposition \ref{WinklerDjiePropo}, let $(u,v)$ be the local solution to system \eqref{sys1}, which blows up at finite time $T_{max}$ in the sense of \eqref{Blow-UpClassicalSense}, and $\Phi(t)$ the energy function defined in Lemma \ref{LemmaFirstStepDerivativePhi}. If, additionally, for some $p_0>\frac{n}{2}(m_2-m_1)$ it is known that for some $L>0$
\begin{equation}
\lVert u(\cdot,t)\rVert_{L^{p_0}(\Omega)}\leq L\quad \textrm{ for all } t \in (0,T_{max}),
\end{equation}
then there exists  $K>0$ with this property: 
\begin{align}\label{BoundednessU^p}
\left\|u(\cdot,t)\right\|_{L^{p}(\Omega)} 
\le K \quad \textrm{ for all } t \in (0,T_{max}),
\end{align}
 and for any $q_1>n+2$  it holds that
\begin{align}\label{BoundednessCrossDiffusion} 
\left\|u(\cdot,t)(u(\cdot,t)+\alpha)^{m_2-2}\nabla v(\cdot,t)\right\|_{L^{q_1}(\Omega)} 
\le K \quad \textrm{ for all } t \in (0,T_{max}).
\end{align}
\begin{proof}
With the results of lemmata  \ref{LemmaFirstStepDerivativePhi}, \ref{LemmaEstimateu^pandu^p+m2-2} and \ref{LemmaBoundednessNormUpCrossDiffusion}  in hour hands, we have by plugging 
\eqref{powerc} into \eqref{ph} 
\begin{align}\label{sum} 
\Phi'(t) 
\le 
-
\left(E_0-E_1\ep \right) 
\int_{ 
 \Omega} \left| \nabla (u+\alpha)^
{\frac{p+m_1-1}{2}}\right|^2
+E_1 E_7+E_5\quad \textrm{ for all } t\in (0,T_{max}).
\end{align}
Moreover, in order to have strictly positivity of the first term  
on the right-hand side  of this gained inequality,  we choose $\ep$ small 
enough as to satisfy $E_0-E_1\ep >0.$ In this way, taking into consideration \eqref{powerb}, relation \eqref{sum} reads 
\begin{align}\label{comp}  
\Phi'(t) 
\le 
-
\left(E_0-E_1\ep \right) 
\left( 
E_6\Phi^\lambda(t) + 1 
\right)
+E_1 E_7+E_5\le 
-J_1 
\Phi^\lambda (t)
+J_2 \quad \textrm{ on } (0,T_{max}),
\end{align}
where $J_1:=\left(E_0-E_1\ep \right) E_6$  and  $J_2:=E_1 E_7+E_5$. Subsequently, we arrive at this initial problem
\begin{equation*}\label{MainInitialProblemWithM}
\begin{cases}
\Phi'(t)\leq J_2-J_1 \Phi^\lambda(t) \quad t \in (0,T_{max}),\\
\Phi(0)=\frac{1}{p}\int_\Omega (u_0+\alpha)^p, 
\end{cases}
\end{equation*}
so to have, by an application of a comparison principle, 
\begin{equation}\label{Boundu^pLemmaWinklerTao}
\Phi(t)\leq \max\left\{\Phi(0),\left(\frac{J_2}{J_1}\right)^\frac{1}{\lambda}\right\}=:L_1\quad \textrm{for all}\quad t\in(0,T_{max}).
\end{equation}
On the other hand, from this bound, elliptic regularity results applied to the second equation of system \eqref{sys1}, i.e. $-\Delta v=u-M$, imply $v\in L^\infty((0,T_{max});W^{2,p}(\Omega))$ and, hence, $\nabla v \in L^\infty((0,T_{max});W^{1,p}(\Omega))$. In particular, the Sobolev embeddings (from Lemma \ref{Lemmapbarra} is $p=\bar{p}>q_1>n+2$) infer $\nabla v\in L^\infty((0,T_{max});L^{\infty}(\Omega))$. Consequently, through the H\"{o}lder inequality with exponents $ q_1(m_2-1)/p$  and $1-q_1(m_2-1)/p$ (again Lemma \ref{Lemmapbarra} ensures that $p>q_1(m_2-1)$), we have on $(0,T_{max})$
\begin{equation*}
\begin{split}
\int_\Omega |u(u+\alpha)^{m_2-2}\nabla v|^{q_1}&
\leq\int_\Omega (u+\alpha)^{q_1(m_2-1)}\lvert \nabla v\rvert^{q_1} \leq   \lVert \nabla v (\cdot, t) \rVert_{L^\infty(\Omega)}^{q_1}|\Omega|^\frac{p-q_1(m_2-1)}{p} \left(\int_\Omega (u+\alpha)^{p}\right)^\frac{q_1(m_2-1)}{p}.
\end{split}
\end{equation*}
Therefore, in view of estimate \eqref{Boundu^pLemmaWinklerTao} we also get 
\[\int_\Omega |u(u+\alpha)^{m_2-2}\nabla v|^{q_1}\leq   \lVert \nabla v (\cdot, t) \rVert_{L^\infty(\Omega)}^{q_1}|\Omega|^\frac{p-q_1(m_2-1)}{p} L_1^\frac{q_1(m_2-1)}{p}=:L_2 \quad \textrm{with}\quad  q_1>n+2, \]
so that \eqref{BoundednessU^p} and \eqref{BoundednessCrossDiffusion} are attained posing $K=\max\{(L_1p)^\frac{1}{p},(L_2)^\frac{1}{q_1}\}$.
\end{proof} 
\end{lem}
\subparagraph*{Proof of Theorem \ref{TheoremMain1}} 
Proposition \ref{WinklerDjiePropo} provides the unique local classical solution $(u,v)$ to system \eqref{sys1} which blows up at finite time $T_{max}>0$. By reduction to the absurd, let $(u,v)$ such that for all $p_0>\frac{n}{2}(m_2-m_1)$ it holds that
\[
\limsup_{t \rightarrow T_{max}} \left\|u(\cdot,t)\right\|_{L^{p_0}
(\Omega)} <\infty;
\] 
then, for  some $L>0$ we get
\[
\left\|u(\cdot,t)\right\|_{L^{p_0}
(\Omega)} \leq L\quad \textrm{ for all } t \in (0,T_{max}).
\]
Now, for for $p=\bar{p}$ given in \eqref{ConstantForTechincalInequality_Barp}, Lemma \ref{LemmaBoundednessNormUpCrossDiffusion} ensures that 
\begin{equation}
\begin{cases}
u\in L^\infty((0,T_{max});L^p(\Omega))& (\textrm{for } p= \bar{p}),\\
u(u+\alpha)^{m_2-2}\nabla v \in L^\infty((0,T_{max});L^{q_1}(\Omega))&\textrm{for all } q_1>n+2. 
\end{cases}
\end{equation}
Hereafter, with the same nomenclature used by Tao and Winkler, $u$ also classically solves in $\Omega \times (0,T_{max})$ problem (A.1) of \cite[Appendix A]{TaoWinkParaPara} for 
\begin{equation*}
D(x,t,u)=(u+\alpha)^{m_1-1},\quad f(x,t)=\chi u(u+\alpha)^{m_2-2}\nabla v,\quad g(x,t)\equiv 0. 
\end{equation*}
In particular, from the boundary condition on $v$, we see that (A.2)--(A.5) and the second inclusion of (A.6) for any choice of $q_2$ are complied.  Moreover, always from the definition of $\bar p$, relations (A.8), (A.9) and (A.10)  of \cite[Lemma A.1.]{TaoWinkParaPara} are also valid, so we have through this lemma that for some $C>0$
\begin{equation*} 
\lVert u(\cdot,t)\rVert _{L^\infty(\Omega)}\leq C \quad \textrm{for all}\quad t \in (0,T_{max}),
\end{equation*}
which is in contradiction to  the fact that the solution $(u,v)$ blows up at finite time $T_{max}.$
\qed 
\section{The ordinary differential inequality for $\Phi(t)$: derivation of lower bounds}\label{SectionODIForLower}
In preparation to the last proof, let us now use some of the above derivations to obtain an ODI for the energy function $\Phi(t) := \frac{1}{\bar{p}} 
		\int_{\Omega} \left(u+\alpha\right)^{\bar{p}}$. This ODI, actually, is satisfied by $\Phi(t)$ both if such energy function is associated to a local or a global solution $(u,v)$ to system \eqref{sys1}; despite this, since we will make use of this ODI to estimate the blow--up time for $T_{max}$, we also confine the forthcoming lemma to the case of unbounded solutions.
\begin{lem}\label{lem diff ineq}
Under the assumptions of Proposition \ref{WinklerDjiePropo}, let $(u,v)$ be the local solution to system \eqref{sys1}, which blows up at finite time $T_{max}$ in the sense of \eqref{Blow-UpClassicalSense},  and $\Phi(t)$ the energy function defined in Lemma \ref{LemmaFirstStepDerivativePhi}. Then there exist $E_8, E_9, E_5>0$ such that $\Phi(t)$ satisfies this ODI
\begin{equation} \label{diff ineq Phi}
\Phi'(t)\leq E_8 \Phi^{\gamma}(t) + E_9\Phi^{\delta}(t) + E_5 \quad \textrm{ on } (0,T_{max}),
\end{equation}
being $\gamma, \delta>1$ as in \eqref{gammabetadelta}.
 \begin{proof}
We start from Lemma \ref{LemmaFirstStepDerivativePhi} and use the Gagliardo-Niremberg inequality to estimate the last term on the right hand side of \eqref{ph}. For $k:=\frac{2(p+m_2-1)}{p+m_1-1}$ as in Lemma \ref{Lemmapbarra} and $a_3$ defined in \eqref{a_1},  if we set   
\[\mathfrak{q}=\mathfrak{s}=\frac{2p}{p+m_1-1}, \mathfrak{p}= k,\] 
Lemma \ref{G-N_InequLemma} yields 
\begin{align}\label{u^p+m_2-1} 
\int_{ \Omega} (u+\alpha)^{p+m_2-1} 
&= \left\| (u+\alpha)^{\frac{p+m_1-1}{2}} \right\|_{L^{k}(\Omega)}^{k}\\\nonumber
&\le
C_{GN}
\left\|
\nabla (u+\alpha)^{\frac{p+m_1-1}{2}}
\right\|_{L^2(\Omega)}^
{ka_3} 
\left\|
(u+\alpha)^{\frac{p+m_1-1}{2}}
\right\|_{L^{\frac{2p}
{p+m_1-1}}(\Omega)}^{
k(1-a_3)}
\\\nonumber
&\quad\, 
+C_{GN}
\left\|
(u+\alpha)^{\frac{p+m_1-1}{2}}
\right\|_{L^
{\frac{2p}{p+m_1-1}
}(\Omega)}^
{k} \quad \textrm{ on } (0,T_{max}).
\end{align} 
Thanks to the first of \eqref{rho<1andgamma>delta>1}, and recalling the definitions in  \eqref{gammabetadelta}, an application of the Young inequality with exponents $\sigma=\frac{ka_3}{2}$ and $1-\sigma$ leads to 
\begin{align}\label{u^p+m_2-1 bis}
\int_{ \Omega} (u+\alpha)^{p+m_2-1} 
&\leq \frac{E_0}{E_1} \int_{\Omega} \left | \nabla(u+\alpha)^{\frac {p+m_1 -1}{2}} \right |^2 + c_4\Phi^{\gamma} + c_5 \Phi^{\delta} \quad \textrm{ for all } t\in (0,T_{max}),
\end{align} 
with  
\begin{equation} \label{c_6, c_7, c_8}
c_4= p^{\gamma}C_{GN}(1-\sigma)\Big({\frac{E_0}{E_1C_{GN}\sigma}\Big)^{-{\frac{\sigma}{1-\sigma}}} }\quad \textrm{and}\quad
c_5= p^\delta C_{GN}.
\end{equation}
Then, by inserting \eqref{u^p+m_2-1 bis} into \eqref{ph} gives the claimed ordinary differential inequality
\begin{align}\label{Phi'}
\Phi'(t) 
&\le 
E_8 \Phi^{\gamma}(t) + E_9\Phi^{\delta}(t) +E_5\quad \textrm{ for all } t\in (0,T_{max}),
\end{align} 
 with $
E_8=c_4 E_1,\ \ E_9=c_5E_1.$
\end{proof}
\end{lem}
\subparagraph*{Proof of Theorem \ref{main}} 
For $n\in \N$ and $m_1,m_2 \in \R$ complying with the blow--up restrictions \eqref{b},  let $p=\bar{p}$ be the number given in Lemma \ref{Lemmapbarra} and $T_{max}$ the finite blow--up time, in  $L^\infty(\Omega)$--norm, of the local solution $(u,v)$ to system \eqref{sys1} provided by Proposition \ref{WinklerDjiePropo}. Since $p=\bar{p}>p_0$, from Theorem \ref{TheoremMain1} we know that $\limsup_{t\rightarrow T_{max}}\frac{1}{p}\int_\Omega (u+\alpha)^p=\infty.$ On the other hand, Lemma \ref{lem diff ineq} ensures that $u$ satisfies the ODI \eqref{diff ineq Phi} for  any $0<t<T_{max}$, where in particular it is seen that the function $\Psi(\xi)=E_8\xi^\gamma+E_9\xi^\delta+E_5$ obeys  the Osgood criterion \eqref{OsgoodCriterion}, where $E_5=E_5(\bar{p}), E_8=E_8(\bar{p}), E_9=E_9(\bar{p})$ have been computed in lemmata \ref{LemmaFirstStepDerivativePhi} and \ref{lem diff ineq} and $\gamma=\gamma(\bar{p})>1,\delta=\delta(\bar{p})>1$ defined in \eqref{rho<1andgamma>delta>1}. Thereafter, by integrating \eqref{diff ineq Phi} between $0$ and $T_{max}$,  we obtain estimate \eqref{blow-up-time}, and the proof is completed. 
\qed
\begin{remark}
We observe that, conversely to what happens with relation \eqref{blow-up-time},  it is possible to obtain an explicit expression for the lower bound $T$ by reducing \eqref{Phi'} as follows: from the definition of $M$, i.e. $M=\frac{1}{|\Omega|}\int_\Omega u_0(x)dx$, and the H\"older inequality we can estimate $E_5$ in relation \eqref{diff ineq Phi} as
\begin{equation} \label{constant}
E_5 =\frac{E_5}{M}\frac{1}{|\Omega|}\int_{\Omega}u_0 \leq \frac{E_5}{M}\frac{1}{|\Omega|} \int_{\Omega} (u+ \alpha) \leq E_{10} \Phi^{\frac 1 p},
\end{equation} 
with 
$$E_{10}=\frac{E_5}{M}\Big(\frac{p}{|\Omega|}\Big)^{\frac 1 p},$$ 
so that \eqref{diff ineq Phi}  can be rewritten in this form:
\begin{equation}\label{ldiff}
\Phi'(t)\leq E_8 \Phi^{\gamma}(t) + E_9 \Phi^{\delta}(t) + E_{10} \Phi(t)^{\frac 1 p}\quad \textrm{ on } (0,T_{max}).
\end{equation} 
Now, similarly to what done in \cite{Marras-Viglia-ComptesGradient}, since $\Phi$ blows up at finite time $T_{max}$ there exists a time $t_1 \in [0, T_{max})$ such that 
\begin{equation*}
\Phi(t) \geq \Phi(0) \quad  \textrm{for all } \, t \geq t_1 \in[0,T_{max}).
\end{equation*}
From  $\gamma >\delta >\frac 1  p $ (recall \eqref{rho<1andgamma>delta>1}), we can estimate the second and third  terms of \eqref{ldiff} by means of $\Phi^{\gamma}$:
\begin{equation} \label{Phi^delta < Phi^gamma}
\Phi^{\delta}(t) \leq  \Phi(0)^{\delta - \gamma} \Phi^{\gamma}(t)\quad \textrm{ and }
\quad \Phi^{\frac 1 p}(t) \leq \Phi(0)^{\frac 1 p - \gamma} \Phi^{\gamma}(t)\quad \textrm{for all } \, t \geq t_1 \in[0,T_{max}).
\end{equation}
By plugging expressions \eqref{Phi^delta < Phi^gamma} into \eqref{ldiff} we obtain for
\begin{equation*}
H=E_8+ E_9 \Phi(0)^{\delta -\gamma} + E_{10} \Phi(0)^{\frac 1 p - \gamma},
\end{equation*}
\begin{equation}\label{ldiff bis}
\Phi'(t)\leq H \Phi^{\gamma}(t) \quad \textrm{for all } \, t \geq t_1 \in[0,T_{max}), 
\end{equation}
so that an integration of \eqref{ldiff bis} on $(t_1,T_{max})$ yields this explicit lower bound for $T_{max}$:
\begin{align*}
\frac{\Phi(0)^{1-\gamma}}{H (\gamma -1)}= \int_{\Phi(0)}^{\infty} \frac{d\tau} {H\tau^{\gamma}}\leq \int_{t_1}^{T_{max}}  d\tau\leq \int_0^{T_{max}} d\tau =  T_{max}. 
\end{align*}
\end{remark}
\section*{Acknowledgments}
The authors would like  to express their sincere gratitude  to Professor Stella Vernier Piro  and Professor Tomomi Yokota  for giving them  the precious opportunity of a joint study and  their encouragement. GV and MM are members of the Gruppo Nazionale per l'Analisi Matematica, la Probabilit\`a e le loro Applicazioni (GNAMPA) of the Istituto Na\-zio\-na\-le di Alta Matematica (INdAM) and are partially supported by the research project \textit{Integro--differential Equations and Non--Local Problems}, funded by Fondazione di Sardegna (2017). 

\end{document}